\documentclass[twoside,12pt]{amsart} 
\usepackage[dvips]{graphicx}
\DeclareGraphicsExtensions{ps,eps}
\makeatletter
\def\serieslogo@{}
\makeatother
\makeatletter 
\def\@setcopyright{}
\makeatother

\newtheorem{Theorem}{Theorem}[section]

\newtheorem{Proposition}[Theorem]{Proposition}

\theoremstyle{definition}

\theoremstyle{remark}

\newcommand{\ETTE}{{{\mathbb{T}}}}
\newcommand{\ERRE}{{{\mathbb{R}}}}
\newcommand{\ZED}{{{\mathbb{Z}}}}

\begin{document}

\title
{
Approximation of a simple Navier-Stokes model by monotonic rearrangement
}

\author{Yann Brenier}
\address{
CNRS, Centre de math\'ematiques Laurent Schwartz,
Ecole Polytechnique, FR-91128 Palaiseau, France.
}
\curraddr{}
\email{brenier@math.polytechnique.fr}

\maketitle
\markboth{}{}

$$ $$
\section*{Abstract}
$$ $$

We consider a very simple Navier-Stokes model for compressible
fluids in one space dimension. 
Neglecting the temperature equation and assuming both the pressure
and the viscosity to be proportional to the density, we get
\begin{equation}
\begin{split}
\partial_t(\rho v)+\partial_x(\rho v^2+\lambda\epsilon \rho))=\epsilon\partial_x(\rho\partial_x v),
\\
\partial_t \rho +\partial_x(\rho v)=0,
\end{split}
\end{equation}
where $\lambda,\epsilon$ are positive constants.
We assume both the density field  $\rho(t,\cdot)>0$ and the velocity field $v(t,x)$ to be 
1-periodic in space, with unit mean for the density.
We denote by $X(t,a)=a+\xi(t,a)$  the "material" positions of the fluid parcels, so that $\xi$ is one-periodic in $a\in \ERRE$,
$\partial_t X(t,a)=v(t,X(t,a))$ and  $\partial_a X(t,a)\rho(t,X(t,a))=1$.
Then, we show that $X(t,a)=a+\xi(t,a)$ can be obtained as the unique limit, as $h>0$ goes
to zero and $nh$ goes to $t$, of the time-discrete approximation
$X_{h,n}(a)=a+\xi_{h,n}(a)$ obtained from the recursion formulae
%%%%%%%%%%%%%%%%%%%%%%%%%%%%%%
\begin{equation}
\begin{split}
X_{h,n}(a)=
[a+(1+h\lambda)\xi_{h,n-1}(a)+hZ_{h,n-1}(a))+\sqrt{2\epsilon h}N(a/h)]]^\sharp
\\
Z_{h,n}(a)=(1-\lambda h)Z_{h,n-1}(a)-\lambda^2\xi_{h,n-1}(a)),
\end{split}
\end{equation}
where $X_{h,0},Z_{h,0}$ are suitably initialized,
$N$ is a 1-periodic function with average $0$ and variance $1$,
and $^\sharp$ denotes the rearrangement  operator 
in increasing order for functions depending on $a\in \ERRE$. Since
the term $N(a/h)$ can be interpreted as a deterministic approximation of
a random variable $W$ with zero mean and unit variance, we see that, in some sense,
the resolution of the Navier-Stokes system has been reduced, 
through the crucial intervention
of the monotonic 
rearrangement operator, to the discretization of a 
trivial stochastic differential equation, namely
$$
d\xi=\sqrt{2\epsilon dt}W+(\lambda \xi+Z)dt,\;\;\;dZ=-(\lambda \xi+Z)\lambda dt.
$$
\\
In addition, our result can be easily extended to a
related Navier-Stokes-Poisson system.

\subsubsection*{}
Key words: fluid mechanics, monotonic rearrangement, optimal transport,
systems of particles
\\
MSC: 35Q35
\newpage

\section{Introduction}

Optimal transport theory has been succesfully used to treat many different types of parabolic 
equations as gradient flows of various ``entropy functionals'' for various 
``transportation metrics'', the canonical example being the regular scalar heat equation viewed
by Jordan, Kinderlehrer and Otto \cite{JKO} 
as the gradient flow of the Boltzmann entropy for the quadratic
Monge-Kantorovich MK2 (frequently nammed Wasserstein) metric. There has been a huge literature on this
subject in the last 20 years (study of the heat equation in a very general framework, 
porous-medium equations, thin-film flow equations,
chemotaxis models, etc.., see, as a sample, \cite{AGS,Vi,BCC,GO}).
Therefore, it is natural to try to treat the Navier-Stokes or the Euler equations of fluid mechanics 
in a similar way. Let us mention some attempts such as 
the introduction of a JKO scheme for the Euler equation by Gangbo
and Westdickenberg \cite{GW} or the related treatment of the Navier-Stokes equation by Gigli and
Mosconi \cite{GMo}. 
However, we are probably still very far to be able to solve the initial value 
problem for the Euler or the Navier-Stokes
equations by optimal transport tools.
Nevertheless, there is more hope for fluid models in one space dimension. Indeed, in one dimension,
optimal transport according to the MK2
metric is rather simple, as well known, 
since it can be entirely rephrased in terms of monotonic rearrangement and 
pseudo-inverse functions.
Very simple models, such as the inviscid Burgers equation (and more generally
multidimensional scalar conservation laws) or one-dimensional pressure-less Euler equations, 
have been successfully written,
using pseudo-inverse functions and monotonic rearrangements,
directly as subdifferential equations with nice contraction properties in all
$L^p$ spaces \cite{BBL,Br7,BGSW}. 
In the present paper, we show that similar ideas can be extended to some
Navier-Stokes models for compressible fluids in one space dimension, leaving widely open the
the much more challenging case of 2D and 3D Navier-Stokes equations.
More precisely,
let us consider the (over-simplified) fluid mechanics model of Navier-Stokes type,
without energy equation, in one space variable:
\begin{equation}
\label{NS}
\begin{split}
\partial_t(\rho v)+\partial_x(\rho v^2+p(\rho))=\partial_x(\mu(\rho)\partial_x v),
\\
\partial_t \rho +\partial_x(\rho v)=0,
\end{split}
\end{equation}
where $\rho=\rho(t,x)>0$, $v=v(t,x)$, 
respectively denote the density and the velocity
of the fluid at time $t>0$ and point $x$, on the real line, and $p$, $\mu$  denote the pressure and the viscosity of the
fluid that are supposed to be just given functions of the density. (The genuine Navier-Stokes
equations would include an additional equation for the temperature field, on which both $p$
and $\mu$ would depend.)
Our approach is restricted to the (rather unphysical) case:
\begin{equation}
\label{isotherm}
\mu(\rho)=\epsilon \rho,\;\;\;p(\rho)=\lambda\epsilon\;\rho.
\end{equation}
[As a matter of fact, if we neglect temperature effects, as we have done, it would be 
physically more consistent to assume the viscosity to be a constant and the pressure to be proportional to the density, but our method would not apply then.]
The $\epsilon$ parameter will be interpreted subsequently as a level of noise in the formulation we
will propose for the Navier-Stokes system.
In particular,
the "pressure-less viscosity-less" case $p=\mu=0$ will correspond to a zero level of noise.
The $\lambda$ parameter just scales the pressure with respect to the viscosity.
\\
\\
Let us now write down the NS model (\ref{NS}), under assumption 
(\ref{isotherm}),
in ``material''(or Lagrangian) coordinates. For that purpose,
we define the material positions $X(t,a)$ of the fluid parcels which satisfy:
\begin{equation}
\label{material}
\partial_t X(t,a)=v(t,X(t,a)),\;\;\;\rho(t,X(t,a))\partial_a X(t,a)=1,\;\;\;X(t,a+1)=X(t,a)+1.
\end{equation}
Then, we find (after standard calculations with 
a multiple use of the chain rule)
\begin{equation}
\begin{split}
\partial_{t} X(t,a)
=V(t,a)-\partial_a(\frac{\epsilon}{\partial_a X(t,a)})
\\
\partial_{t} V(t,a)
=-\partial_a(\frac{\lambda\epsilon}{\partial_a X(t,a)})
\end{split}
\end{equation}
or, equivalently,
\begin{equation}
\label{NS material0}
\begin{split}
X(t,a)=a+\xi(t,a),
\\
\partial_{t} \xi(t,a)+\partial_a(\frac{\epsilon}{\partial_a X(t,a)})
=Z(t,a)+\lambda \xi(t,a)
\\
\partial_{t} Z(t,a)
=-\lambda Z(t,a)
-\lambda^2 \xi(t,a)
\end{split}
\end{equation}
by substituting $Z$ for $V-\lambda\xi$.
Our key observation is that this system looks like a very mild modification of
\begin{equation}
\begin{split}
X(t,a)=a+\xi(t,a),
\\
\partial_{t} \xi(t,a)+\partial_a(\frac{\epsilon}{\partial_a X(t,a)})
=0,
\end{split}
\end{equation}
which is nothing but than the ``material'' version of the linear heat equation
(written in ``JKO style'' \cite{JKO})
\begin{equation}
\partial_t \rho +\partial_x(\rho v)=0,\;\;\;v=-\frac{\partial_x \rho}{\rho},
\end{equation}
through the same change of variables (\ref{material}).
This suggests the following discrete scheme, based on the simulation of the
heat equation by random walk, to approximate the 
Navier-Sokes equation
written in material coordinates (\ref{NS material0}).
\\
Given a uniform time step $h>0$ and a uniform grid on $\ERRE$ with mesh size $1/M$.
we approximate
$X(t,a)$, at each discrete time $t=nh$, $n=0,1,2,\cdot\cdot\cdot$,
by a non-decreasing sequence 
$k\rightarrow X_{n,k}$ corresponding to each subinterval $(k-1)/M<a<k/M$, 
for $k\in\ZED$. Because of spatial periodicity, it is convenient to introduce
\begin{equation}
\label{xi}
\xi(t,a)=X(t,a)-a,
\end{equation}
which is 1-periodic in $a$.
Accordingly, we require $X_{n,k+M}=1+X_{n,k}$ (which makes
the scheme effective for computation) and introduce $\xi_{n,k}=X_{n,k}-k/M$.
We first perform a predictor step 
\begin{equation}
\label{fully}
\begin{split}
\hat X_{n+1,k}=\hat \xi_{n+1,k}+k/M
\\
\hat \xi_{n+1,k}=
(1+h\lambda)\xi_{n,k}+hZ_{n,k}+\sqrt{2\epsilon h}N_{n,k}
\\
Z_{n+1,k}=(1-\lambda h)Z_{n,k}-h\lambda^2\xi_{n,k},
\end{split}
\end{equation}
where $N_{n,k}$ is a sequence, $M$-periodic in $k$, of $M$ independent random numbers,
with  expectation $0$ and unit variance.
So far, this can be seen as a discretization of the rather trivial 
stochastic differential equation
\\
\begin{equation}
\label{stocha}
\begin{split}
d\xi=(Z+\lambda \xi)dt+\sqrt{2\epsilon dt}\; W
\\
\frac{dZ}{dt}+\lambda Z=-\lambda^2 \xi,
\end{split}
\end{equation}
where
$W$ is a normalized  white noise.
However, we observe that there is no reason 
for  $k\rightarrow \hat X_{n+1,k}$ to be increasing.
So, we perform a corrector step just by sorting this sequence in increasing order,
which is enough to define the next sequence  $k\rightarrow  X_{n+1,k}$. As discussed below,
we even simplify the scheme by using a deterministic simulation of the white noise, for instance
by substituting for the random number $N_{n,k}$ the value $+1$ if $k$ is even, 
and $-1$ if $k$ is odd, independently of $n$. This simplification is possible, because
the corrector step  generates a large amount of mixing.
Let us point out that the sorting step is crucial to lift the rather uninteresting
stochastic ordinary differential system (\ref{stocha}) to the level of the more intriguing
Navier-Stokes equations!
\\
It is now very easy to define the natural continuous version of this scheme with respect to the space variable $a$, while keeping
$t$ discrete. To do that we use the concept of rearrangement of functions in increasing order, which
is just the continuous version of the sorting operator for sequences of
real numbers. However, we have to be a little bit careful to handle the required space periodicity.
Let us consider a given
function $Y \in I+L^2(\ETTE)$, where $I$ denotes the identity map on $\ERRE$ and $L^2(\ETTE)$ is the
space of all locally square integrable 1-periodic functions on $\ERRE$. 
Then, its rearrangement in increasing order is the unique function
$Y^\sharp\in I+L^2(\ETTE)$ such that
\begin{equation}
\label{operator}
\begin{split}
\partial_a Y^\sharp\ge 0\;,\;\;\;\;
\int_0^1 G(Y^\sharp(a))da=
\int_0^1 G(Y(a)da,
\end{split}
\end{equation}
for all 1-periodic continuous function $G$, $G\in C^0(\ETTE)$.
Then, the semi-discrete version of the scheme 
just reads (with obvious notation): 
\begin{equation}
\label{scheme}
\begin{split}
\hat X_{n+1}(a)=\hat \xi_{n+1}(a)+a,
\\
\hat \xi_{n+1}(a)=
(1+h\lambda)\xi_{n}(a)+hZ_{n}(a)+\sqrt{2\epsilon h}W
\\
Z_{n+1}(a)=(1-\lambda h)Z_{n}(a)-h\lambda^2\xi_{n}(a),
\\
X_{n+1}=[\hat X_{n+1}]^\sharp,
\end{split}
\end{equation}
where $W$ is a normalized white noise.
As a matter of fact, we rather use a deterministic approximation of the white noise term, namely
\begin{equation}
\label{pseudo}
W=N(a/r),
\end{equation}
where $N$ is a fixed 1-periodic bounded function with zero mean and unit variance and $1/r$ is
a sufficiently large integer.
This prevents us
to rely on any stochastic analysis. 
\\
Let us point out, before stating our convergence result, that, due to its very definition,
the semi-discrete scheme immediately becomes a fully discrete scheme, under the following conditions.
We first approximate the initial
data $(X_0(a),Z_0(a))$ and the identity map $a$
by piecewise constant functions in $a$,
on a uniform grid, with constant values in each interval $](k-1)/M,k/M[$, for $k\in\ZED$.
Observe that the uniformity of the grid is crucial, 
because the rearrangement operator in increasing order preserves
piecewise constant functions only on a uniform grid.
We also crucially need the
"noise function" $N(a/r)$  to be itself piecewise constant on the $same$ grid. Let us consider for
instance the simplest noise function $N(a)=-1$ for $0<a<1/2$ and $N(a)=1$ for $1/2<a<1$.
Since $L=1/r$ is supposed to be an integer, we see that $N(a/r)$ will be piecewise constant on the
grid of mesh width $1/M$ as soon as we choose $M$ to be a multiple of $2L$. If we choose, with
economy, $M=2L$, we see that $N(a/r)$ will take value $-1$ on each interval $](2k-2)/M,(2k-1)/M[$
and $+1$ on each interval $](2k-1)/M,2k/M[$, for $k\in\ZED$.
Then the time-discrete scheme, without modification, coincides with the fully-discrete scheme discrete 
(\ref{fully}) where we substitute for the random number $N_{n,k}$ the value $+1$ if $k$ is even, 
and $-1$ if $k$ is odd.
Let us also observe that, by construction, the fully discrete scheme can be trivially coded
(with the help of a fast sorting algorithm, taking into account periodicity), and the resulting complexity, to reach the solution at a finite fixed time $t$, is of order $O(h^{-1}r^{-1}\log(1/r))$ on a sequential machine.
\\
It is now time to state our main result.

\begin{Theorem}
\label{main}
Let us consider the time-discrete scheme (\ref{scheme},\ref{pseudo}), 
where $N$
is the 1-periodic binary function with value $-1$ for $0<a<1/2$ and
$+1$ for $1/2<a<1$, and $1/r$ an integer large enough so that $r=O(h)$. 
Let us consider initial conditions $X_0=I+\xi_0$, $Z_0$, where $I$ is the identity map
and $\xi_0$ and $Z_0$ are smooth 1-periodic functions, with $\inf X_0'(a)>0$.
Let us linearly interpolate in time
the $(X_{h,n},Z_{h,n})$ and denote the result by $(X_h,Z_h)$. Then, as $h\rightarrow 0$,
the entire family converge to  a solution $(X,Z)$ of the material formulation of
the Navier-Stokes system (\ref{NS material0}), under assumption (\ref{isotherm}).
\end{Theorem}

Symbolically, to describe this limit, we suggest the following continuous formulation:
\begin{equation}
\label{new}
\begin{split}
X(t+dt,a)=[X(t,a)+(Z(t,a)+\lambda (X(t,a)-a))dt+ \sqrt{2\epsilon dt}\; W]^\sharp
\\
\frac{dZ}{dt}+\lambda Z=-\lambda^2 (X-a),
\end{split}
\end{equation}
where $W$ stands for a normalized white noise
and $^\sharp$ denotes the rearrangement  operator
in increasing order.
\\
\\
%ZZZZZZZZZZZZ
Our method of proof will involve three main steps:
\\
i) The Navier-Stokes system, written in material coordinates, is first interpreted as a modified gradient flow with delay terms for a suitable convex functional, in a slightly more general situation. As a byproduct, we establish  the
global $L^2$ well-posedness of the system, in a general setting of non smooth initial conditions
(including vacuum, i.e. zones where the density field vanishes).
\\
ii) In the special case (\ref{isotherm}),
we show that the previous system can be further reduced to a non-autonomous viscous scalar
conservation laws, for which existence of global smooth solutions can be routinely established.
\\
iii) Finally, we observe that our time-discrete scheme is an adaptation of the "transport-collapse-kinetic" method \cite{Br0,Br1,GMi,Br2,Br4} to the discretization of a suitable coupled scalar parabolic equation (\ref{scalar3})
and prove the desired  convergence theorem.
%ZZZZZZZZZZZ
\\
\\
Let us conclude this introduction by mentioning that our approach can be adapted to some interesting cases when a self-consistent forcing term is added to the NS equations.
This is in particular true for the Navier-Stokes-Poisson (NSP) system

\begin{equation}
\label{NSP}
\begin{split}
\partial_t(\rho v)+\partial_x(\rho v^2+p(\rho))=
\epsilon\partial_x(\mu(\rho)\partial_x v)
-\rho\partial_x\phi,
\\
\partial_t \rho +\partial_x(\rho v)=0,\;\;\;
\beta\partial^2_{xx}\phi=\rho-1,
\end{split}
\end{equation}

where $\beta$ is a constant, $\phi=\phi(t,x)$ can be interpreted ad a gravitational potential if 
$\beta>0$, or as
an electrostatic potential if $\beta<0$, with a "neutralizing"
background of unit mean.
In material coordinates,  this system reads
\begin{equation}
\label{NSP material}
\partial^2_{tt} X(t,a)=\partial_a [-p(1/\partial_a X)
+\nu(\partial_a X)\partial^2_{at}X(t,a)]
+F(X(t,a)-a),
\end{equation}
where 
\begin{equation}
\label{F NSP}
F(y)=\frac{y}{\beta}.
\end{equation}
For this Navier-Stokes-Poisson system, under assumption (\ref{isotherm}),
we just have to modify the noisy differential system in a straightforward way
\begin{equation}
\label{stocha NSP}
\begin{split}
d\xi=(Z+\lambda \xi)dt+\sqrt{2\epsilon dt}\; W
\\
\frac{dZ}{dt}+\lambda Z=F(\xi)-\lambda^2 \xi,
\end{split}
\end{equation}
where  $F$ is given by (\ref{F NSP}).
So we may immediately include the Navier-Stokes-Poisson case, just by adding a smooth forcing term,
$F(y)$ with bounded derivative in $y$ and $F(0)=0$, to the Navier-Stokes equations, as in (\ref{NSP material}), 
without necessarily assuming that $F$ has form (\ref{F NSP}). This is what we will do through the paper.

\section{Reduction to a gradient flow with delay}

In this section, we consider the Navier-Stokes model
(written in material coordinates) (\ref{NSP material}),
in the case
when the viscosity coefficient $\mu$ and
the pressure $p$ are linked to each other through a given smooth strictly 
convex function $\tau>0\rightarrow \theta(\tau)$, in the following way:
\begin{equation}
\label{etat}
\mu(1/\tau)=\epsilon\; \tau \theta"(\tau)
\;\;\;
p(1/\tau)=-\lambda\epsilon\; \theta'(\tau),
\;\;\;\forall \tau>0.
\end{equation}
for some constants $\epsilon,\lambda>0$. 
We recover the case (\ref{isotherm}) discussed in the introduction
for the special choice:
\begin{equation}
\label{theta}
\theta(\tau)=-\log \tau,\;\;\theta'(\tau)=-1/\tau,\;\;\theta"(\tau)=1/\tau^2.
\end{equation}
Another interesting example is $\theta(\tau)=\tau(\log\tau-1)$,
for which $p=-\lambda\epsilon\; \log\tau$, and
$\mu=\epsilon$ (i.e. constant viscosity).
Under  assumption (\ref{etat}), the NS equation written in material coordinates, with forcing $F$,
(\ref{NSP material}) becomes
\begin{equation}
\partial^2_{tt} X(t,a)
=\lambda\epsilon\;\partial_a(\theta'(\partial_a X(t,a)))
+\epsilon\;\partial_a(\theta"(\partial_a X(t,a))\partial^2_{ta}X(t,a))+F(X(t,a)-a)
\end{equation}
which can be immediately blown up as a first order system:
\begin{equation}
\begin{split}
\partial_{t} X(t,a)
=V(t,a)+\epsilon\;\partial_a(\theta'(\partial_a X(t,a)))
\\
\partial_{t} V(t,a)
=\lambda\epsilon\;\partial_a(\theta'(\partial_a X(t,a)))+F(X(t,a)-a)
\end{split}
\end{equation}
or, equivalently
\begin{equation}
\label{NS theta}
\begin{split}
X(t,a)=a+\xi(t,a),
\\
\partial_{t} \xi(t,a)
=Z(t,a)+\lambda \xi(t,a)+\epsilon\;\partial_a(\theta'(\partial_a X(t,a)))
\\
\partial_{t} Z(t,a)
=-\lambda Z(t,a)
-\lambda^2 \xi(t,a)+F(\xi(t,a))
\end{split}
\end{equation}
just by introducing $Z=V-\lambda (X-a)$, i.e. $V=Z+\lambda (X-a)$.
The resulting system is a mild modification of a gradient flow for a convex functional.
More precisely, let us introduced the Hilbert space $H=L^2(\ETTE)$ and the closed convex
cone $K$  of all non-decreasing functions in $I+H$. 
\\
Now let us introduce the lsc convex functional
$\xi \in H\rightarrow \Theta[\xi]\in ]-\infty,+\infty]$, with value  $+\infty$ if $X=I+\xi$ is not in $K$ and 
\begin{equation}
\label{Theta}
\Theta[\xi]=\int_0^1\theta(1+\xi'(a))da\in ]-\infty,+\infty],
\end{equation}
otherwise. 
[Notice that whenever $X=I+\xi \in K$, its derivative $X'(a)=1+\xi'(a)$ can be seen as a locally bounded nonnegative Borel measure on $\ETTE$. Then $\Theta[\xi]$ can be
 more precisely defined as "a convex function of measures" by Legendre duality, as done, for example,
 in  \cite{DT}:
\begin{equation}
\label{Thetabis}
\Theta[\xi]=\sup_f \int_{\ETTE}\{f(a)d(a+\xi(a))-\theta^*(f(a))da\}, 
\end{equation}
where $f$ describes the set of continuous function on $\ETTE$ and $\theta^*$ is the Legendre-Fenchel transform of $\theta$: 
$$
\theta^*(s)=\sup_{\tau\in\ERRE}(s\tau-\theta(\tau)) .
$$
For example, in the case (\ref{theta}) we are mostly interested in, $\theta(\tau)=-\log\tau$ (extended by $+\infty$ for $\tau\le 0$), 
we find $\theta^*(s)=-1-\theta(-s)=-1-\log(-s)$ and we can rewrite
(\ref{Thetabis}) as
$$
\Theta[X]=\sup_g
\int_{\ETTE}\{-\exp(g(a))d(a+\xi(a))+(g(a)+1)da\},
$$
for all continuous functions $g$ on $\ETTE$, by setting $f(a)=-\exp(g(a))$.]
\\
So, the first order system (\ref{NS theta}) can be viewed, in the Hilbert space $H$, as a sub-differential inclusion coupled to a linear ODE in $Z$, namely:
\begin{equation}
\label{NS subdiff}
\begin{split}
-\frac{d\xi}{dt}+Z+\lambda \xi\in\;\; \epsilon\partial \Theta[\xi],
\\
\frac{dZ}{dt}+\lambda Z=-\lambda^2 \xi+F(\xi).
\end{split}
\end{equation}
By  integrating out $Z$ from the second equation,
we may view the system as a gradient flow with 
delay terms (which are harmless since we have assumed $F(y)$ to have bounded derivative in 
$y$ and $F(0)=0$, as in the Navier-Stokes-Poisson case (\ref{F NSP})):

\begin{equation}
\label{NS subdiff2}
-\frac{d\xi}{dt}+Z(0)e^{-\lambda t}
+\int_0^t e^{-\lambda (t-t')}(F(\xi(t'),a)-\lambda^2\xi(t'))dt'
+\lambda \xi \in\;\; \epsilon\partial \Theta[\xi].
\end{equation}

Such formulations are very useful to define global "generalized" solutions for a large set of
initial conditions, with well-posedness in $L^2$, just by using, in a routine way, the general theory of maximal monotone
operators \cite{Brz}. In addition, in this framework, the limit $\epsilon\rightarrow 0$ is trouble-free.
Indeed, in the sense of maximal monotone operator theory, the limit of $\epsilon\Theta[\xi]$ is just the indicator function $1_K[X]$ of the cone $K$, with value $0$ whenever $X=I+\xi$ belongs to $K$ and
$+\infty$ otherwise. Thus, we may conclude this section by asserting:

\begin{Theorem}
\label{subdiff}
Let $(\xi_0,Z_0)$ given in the Hilbert space $H=L^2(\ETTE)$, with $X_0=I+\xi_0$ in the closed convex cone $K$ of all non-decreasing functions in $I+H$.
Then, there is a unique global solution $t \rightarrow (\xi(t),Z(t))\in H\times H$ to (\ref{NS subdiff}),
depending continuously in $H$ on both $t$ and the initial condition
$(\xi(0)=\xi_0,Z(0)=Z_0)$.
This solution can be viewed as a "generalized" solution to the Navier-Stokes
system (\ref{NSP material}),  written in material coordinates, with pressure and viscosity laws given by (\ref{etat}),
and, in particular for (\ref{isotherm}).
In addition, as $\epsilon\rightarrow 0$, these solutions converge to those of:
\begin{equation}
\label{Euler  subdiff}
\begin{split}
-\frac{d\xi}{dt}+Z+\lambda \xi\in\;\;\;\partial 1_K[I+\xi],
\\
\frac{dZ}{dt}+\lambda Z=-\lambda^2 X+F(\xi).
\end{split}
\end{equation}
where $1_K[X]$ is the indicator function of the cone $K$, 
with value $0$ whenever $X=I+\xi$ belongs to $K$ and $+\infty$ otherwise. 
\end{Theorem}

\subsection*{Remark 1}
Formulations of type (\ref{Euler subdiff}) have already been introduced in the absence of viscosity for the pressure-less Euler equations and some of their variants. Let us quote the model of order-preserving vibrating strings, the Chaplygin gas, the pressure-less Euler-Poisson system and also multi-dimensional scalar conservation laws etc ... (\cite{Br5},\cite{Br7},\cite{NS},\cite{Br8}, for example).
In all of these earlier works, $\lambda$ has never been taken into account. It is unclear to us,
whether or not $\lambda$ plays (or should play) a role in these zero viscosity limit equations. We leave that for further investigations.

\subsection*{Remark 2}
Subdifferential formulation (\ref{NS subdiff}), or (\ref{NS subdiff2}),
leads to stability of solutions with respect to their initial condition not
only in $L^2$ but also in all $L^q$ spaces. Indeed, let us consider, on a fixed finite time interval $[0,T]$, 
two solutions $(\xi,Z)$ and $(\tilde \xi,\tilde Z)$. Using formulation 
(\ref{NS subdiff2}), we get
$$
e'(t)\le c\{ ||Z(0)-\tilde Z(0)||^q_{L^q}+e(t)+\int_0^t e(t')dt'\},\;\;\;
e(t)=||\xi(t)-\tilde \xi(t)||^q_{L^q},
$$
where $c$ depends only on $q,T,\lambda$ and the Lipschitz constant of $F$, but not on $\epsilon$.
This easily leads to the $L^q$ stability estimate, uniformly on $\epsilon$,
\begin{equation}
\label{stability Lp}
\sup_{t\in [0,T]}||(\xi,Z)(t)-(\tilde \xi,\tilde Z)(t)||_{L^q}
\le c||(\xi,Z)(0)-(\tilde \xi,\tilde Z)(0)||_{L^q},
\end{equation}
with (another) $c$ depending only on $q,T,\lambda,F$,
but not on $\epsilon$.

\section{Reduction to a scalar parabolic equation with mild coupling}

The sub-differential formulation (\ref{NS subdiff}) of the Navier-Stokes equations in material coordinates
(\ref{NSP material}), that we have introduced in the previous section,
is very powerful since it provides $L^2$ (and, more generally $L^q$) well-posedness for very general data (including vacuum).  However, it is
not clear that the corresponding solutions are smooth, for appropriate initial conditions. 
At least in the special case (\ref{theta}) where $\theta(\tau)=-\log\tau$, this can be done by
reducing (\ref{NSP material}) to the scalar parabolic equation
\begin{equation}
\label{scalar2}
\partial_t u+(Z(t,u)+(x-u)\lambda)\partial_x u=\epsilon\partial^2_{xx} u,
\end{equation}
where $Z$ is coupled to $u$ in a suitable way.
For this purpose, let us introduce the inverse function $u$ of $X=I+\xi$
\begin{equation}
\label{u}
u(t,X(t,a))=a.
\end{equation}
At the moment, let us assume $X(t,a)=a+\xi(t,a)$ to be smooth with $\partial_a X>0$ (which will be justified a posteriori in the case (\ref{isotherm})). We have:
\begin{equation}
\label{prop u}
\begin{split}
u(t,X(t,a))=a,\;\;\;(\partial_x u)(t,X(t,a))\partial_a X(t,a)=1,\;\;
\\
(\partial_t u)(t,X(t,a))=-\partial_x u(t,X(t,a))\partial_t X(t,a).
\end{split}
\end{equation}

Let us now point-wise multiply  the first equation of  (\ref{NS theta}) by
$$
(\partial_x u)(t,X(t,a))=1/\partial_a X(t,a).
$$
Using (\ref{prop u}), we get
$$
-(\partial_t u)(t,X(t,a))=[Z(t,u(t,X(t,a))+\lambda (X(t,a)-a)](\partial_x u)(t,X(t,a))
$$
$$
+\epsilon \partial_a [\theta'(1/(\partial_x u)(t,X(t,a))]/\partial_a X(t,a),
$$
and therefore, by substituting $x$ for $X(t,a)$:
$$
-\partial_t u(t,x)=[Z(t,u(t,x))+\lambda (x-u(t,x)]\partial_x u(t,x)
+\epsilon \partial_x [\theta'(1/\partial_x u(t,x))].
$$
Thus, we get for $(u,Z)$ the following system:
\begin{equation}
\label{scalar1}
\begin{split}
\partial_t u(t,x)+[Z(t,u(t,x))+\lambda (x-u(t,x))]\partial_x u(t,x)=
-\epsilon \partial_x [\theta'(1/\partial_x u(t,x))].
\\
\partial_t Z(t,a)+\lambda Z(t,a)=F(X-a)-\lambda^2 (X(t,a)-a),\;\;\;u(t,X(t,a))=a..
\end{split}
\end{equation}
By integrating out $Z$, we get, more explicitly,
\begin{equation}
\label{flux}
Z(t,a)=Z(0,a)e^{-\lambda t}
+\int_0^t e^{-\lambda (t-t')}(F(X(t',a)-a)-\lambda^2(X(t',a)-a))dt'.
\end{equation}
Next, under assumption (\ref{etat},\ref{theta}),  i.e. $\theta(\tau)=-\log\tau$, we find
\begin{equation}
\label{scalar3}
\begin{split}
\partial_t u(t,x)+[Z(t,u(t,x))+\lambda (x-u(t,x))]\partial_x u(t,x)=\epsilon\partial^2_{xx} u,
\\
\partial_t Z(t,a)+\lambda Z(t,a)=F(X-a)-\lambda^2 (X(t,a)-a),\;\;\;u(t,X(t,a))=a.
\end{split}
\end{equation}
Observe that the coupling of $Z$ with $u$, through $X$ is very mild, 
so that the theory for this system
differs very little from the standard theory of viscous conservation laws, such as
\begin{equation}
\label{scalar0}
\partial_t u+\partial_x(f(u))=\epsilon\partial^2_{xx} u.
\end{equation}
So, without entering in details,  we claim that
for $u_0$ given in the class $C$ of all orientation-preserving 
diffeomorphisms of $\ERRE/\ZED$ and $Z_0$ smooth 1-peridoic function,
then there is a unique smooth solution $(u(t,\cdot)\in C,,Z(t,\cdot))$ to 
the coupled scalar parabolic equation (\ref{scalar3}).

\subsection*{$L^1$ stability with respect to initial conditions}
In addition, thanks to the 
the subdifferential formulation (\ref{NS subdiff}),
we get an $L^1$ stability result for every pair 
$(u,Z)$ and $(\tilde u,\tilde Z)$ of solutions to 
the coupled scalar parabolic equation (\ref{scalar3}).
on every fixed finite time interval $[0,T]$:
\begin{equation}
\label{stability L1}
\sup_{t\in [0,T]}(||u(t)-\tilde u(t)||_{L_1}+||Z(t)-\tilde Z(t)||_{L_1})
\le c(||u(0)-\tilde u(0)||_{L_1}+||Z(0)-\tilde Z(0)||_{L_1})
\end{equation}
with $c$ depending only on $T,\lambda,F$,
but not on $\epsilon$. This immediately follows from 
the $L^q$ stability estimate (\ref{stability Lp}) for the subdifferential formulation (\ref{NS subdiff}),
in the special case $q=1$.

\section{Approximation by a "transport-collapse-kinetic" method with noise}

The transport-collapse-kinetic method is a time-discrete scheme for hyperbolic and viscous conservation
laws introduced and studied in \cite{Br0,GMi,Br2,Br4}.  It was suggested to the author by 
an earlier algorithm proposed by A. Chorin \cite{Ch} for reaction-diffusion equations. 
This scheme turns out to be very well suited for the formulations 
(\ref{scalar3}) and (\ref{NS subdiff}) we have obtained  in the previous sections for
the Navier-Stokes model (\ref{NSP material}). 
\\
 Let us now describe the time-discrete scheme:
\\
Given a uniform time step $h>0$, we denote by
$(X_{h,n}(a)=a+\xi_{h,n}(a),Z_{h,n}(a))$ the approximate solution
at time $t=nh$, for $n=0,1,2,\cdot\cdot\cdot$,
defined in two steps as follows.
\\
We first define $Z_{h,n+1}$ and $\hat X_{h,n+1}=a+\hat \xi_{h,n+1}(a)$ by
\begin{equation}
\label{TCM1}
\begin{split}
\hat \xi_{h,n+1}(a)=
(1+h\lambda)\xi_{h,n}(a)+hZ_{h,n}(a))+\sqrt{2\epsilon h}N(a/r)
\\
Z_{h,n+1}(a)=(1-\lambda h)Z_{h,n}(a)+h(F(\xi_{h,n}(a))-\lambda^2\xi_{h,n}(a)),
\end{split}
\end{equation}
where $N$ is a 1-periodic function with average $0$ and variance $1$,
and $1/r$ is a large integer that may depend on $h>0$. 
[Let us recall that this first step has been designed to approximate the noisy ordinary differential
system (\ref{stocha}), namely:
\begin{equation}
\begin{split}
d\xi=(Z+\lambda \xi)dt+\sqrt{2\epsilon dt}\; W
\\
\frac{dZ}{dt}+\lambda Z=F(\xi)-\lambda^2 \xi.
\end{split}
\end{equation}
where $W$ is a normalized white noise,
with a deterministic approximation to $W$.]
 \\
Next, we crucially introduce a rearrangement step.
Namely, we rearrange $\hat X_{h,n+1}(a)=a+\hat \xi_{h,n+1}(a)$ in increasing order with respect
to $a$, and obtain $X_{h,n+1}(a)=a+\xi_{h,n+1}(a)$.
In other words, we set
\begin{equation}
\label{TCM2}
X_{h,n+1}=\hat X_{h,n+1}^\sharp,
\end{equation}
where we denote by $Y\rightarrow Y^\sharp$ the rearrangement in increasing order, defined
by (\ref{operator}), for functions differing from the identity map $I$ by a $L^q$ 1-periodic function
(i.e. belonging to $I+L^q(\ETTE)$).
Of course, the resulting scheme contains the scheme (\ref{scheme}) discussed in the
introduction, up to the addition of the (harmless) term $F$, that we suppose, as previously
mentioned, smooth with bounded derivative.
Notice that the scheme is almost translation invariant in the variable $a$. Indeed, the rearrangement step is clearly translation invariant, and, in the predictor step, the only non-autonomous term is the "noise" term which is proportional to $N(a/r)$. Thus the scheme is invariant by
any translation $a\rightarrow a+k/r$, with $k\in\ZED$. This observation will be very important for our
analysis. 
We also  observe that $Z_h$ is  entirely "slaved" by $\xi_h$ and depends linearly on it through 
(\ref{TCM1}), so we are only concerned about monitoring $\xi_h$. (Of course we use that 
$F$ is smooth and Lipschitz with $F(0)=0$.)
For the subsequent analysis, it is also simpler
to restrict ourself to the case when  $N$ is the binary function with value $-1$ for $0<a<1/2$ and
$+1$ for $1/2<a<1$. Indeed, in that case, remarkably enough, $(\xi_0,Z_0)=(0,0)$, is a fixed point of the scheme, as $F=0$. [This can be easily checked graphically: we see the exact compensation between the addition of the "noise" term $\sqrt{2\epsilon h}N(a/r)$ and the rearrangement step.]

Let us now prove the following convergence result (from which our main result Theorem \ref{main} immediately follows as the particular case when the force term $F$ is absent):

\begin{Theorem}
Let us consider the time-discrete (\ref{TCM1},\ref{TCM2}),  
where $N$
is the 1-periodic binary function with value $-1$ for $0<a<1/2$ and
$+1$ for $1/2<a<1$ and $1/r$ an integer large enough so that $r=O(h)$. 
Let us consider initial conditions $X_0=I+\xi_0$, $Z_0$, where $I$ is the identity map
and $\xi_0$ and $Z_0$ are smooth 1-periodic functions, with $\inf X_0'(a)>0$.
Let us linearly interpolate in time
the output of the semi-discrete scheme $(X_{h,n},Z_{h,n})$ and denote the result by $(X_h,Z_h)$. Then, as $h\rightarrow 0$,
the entire family converge to $(X,Z)$ solution of 
the material formulation of
the Navier-Stokes system (\ref{NSP material}), under assumption (\ref{isotherm}).
\end{Theorem}

\subsection*{Proof}

\subsubsection*{Sup-norm and Lipschitz estimates}

The rearrangement operator is a non-expansive map on all spaces $I+L^q(\ETTE)$, 
in particular for $q=\infty$. So we can compare in sup-norm any fixed solution $(\xi_{h,n},Z_{h,n})(a)$ 
of the time-discrete scheme to:
\\
i) $(\xi_{h,n},Z_{h,n})(a+kr)$, its space translation by $kr$, for any integer $k\in Z$, which is also solution of the scheme (due to translation invariance, as already mentioned);
\\
ii) the fixed point solution already  $(\xi,Z)=(0,0)$, obtained in the case 
$F=0$,
Let us perform the first comparison. We easily get:
\begin{equation}
\label{translation}
\begin{split}
||(\xi_{h,n},Z_{h,n})-(\xi_{h,n},Z_{h,n})(\cdot+kr)||_\infty
\le (1+hc)^n||(\xi_{0},Z_{0})-(\xi_{0},Z_{0})(\cdot+kr)||_\infty
\\
\le |k|r Lip(\xi_{0},Z_{0})\exp(nhc),
\end{split}
\end{equation}
where Lip denotes Lipschitz constants and
$c$ is a constant depending only on  $\lambda$ and the Lipschitz constant of $F$.
[Indeed, we first get, for the predictor step, 
$$
||(\hat \xi_{h,n+1},Z_{h,n+1})-(\hat \xi_{h,n+1},Z_{h,n+1})(\cdot+kr)||_\infty
$$
$$
\le (1+ch)||(\hat \xi_{h,n+1},Z_{h,n+1})-(\hat \xi_{h,n+1},Z_{h,n+1})(\cdot+kr)||_\infty,
$$
where $\epsilon$ and the "noise" term 
do not play any role, since the translation is an integer multiple of $r$.
Next, using the contraction property of the rearrangement operator, we get
$$
||\xi_{h,n+1}-\xi_{h,n+1}(\cdot+kr)||_\infty
\le ||\hat \xi_{h,n+1}-\hat \xi_{h,n+1}(\cdot+kr)||_\infty.
$$
which is enough to get (\ref{translation}.]
As a matter of fact, we can improve  estimate (\ref{translation})
using the fact that $X_{h,n}(a)=a+\xi_{h,n}(a)$
is increasing in $a$.
Indeed, we get, for all $\omega\in [kr-r,kr]$, 
$$
\xi_{h,n}(a+\omega)-\xi_{h,n}(a)=X_{h,n}(a+\omega)-X_{h,n}(a)
\le X_{h,n}(a+kr)-X_{h,n}(a).
$$
Thus,
$$
||\xi_{h,n}-(\xi_{h,n}(\cdot+\omega)||_\infty
\le  |k|r Lip(\xi_{0},Z_{0})\exp(nhc)\le (|\omega|+r) Lip(\xi_{0},Z_{0})\exp(nhc).
$$
Since $Z_{h,n}$ is slaved by $X_{n,h}$, we have obtained
\begin{equation}
\label{lipschitz}
\begin{split}
||(\xi_{h,n},Z_{h,n})-(\xi_{h,n},Z_{h,n})(\cdot+\omega)||_\infty
\le (|\omega|+r) Lip(\xi_{0},Z_{0})\exp(nhc),
\end{split}
\end{equation}
for all $\omega\in\ERRE$, where
$c$ is a constant depending only on  $\lambda$ and the Lipschitz constant of $F$
and not on $\epsilon$.
This is our first key estimate for the scheme.
Next, let us compare the solution to the fixed point of the scheme
$(\xi=0,Z=0)$, already mentioned and obtained for $F=0$. For this fixed point,
the predictor step provides the values $(\sqrt{2\epsilon h}N(a/r),0)$
Thus, we get
$$
||(\xi_{h,n+1},Z_{h,n+1})-(0,0)||_\infty
\le ||(\hat \xi_{h,n+1},Z_{h,n+1})-(\sqrt{2\epsilon h}N(\cdot/r),0)||_\infty
$$
(because the rearrangement operator is non-expansive)
$$
\le (1+ch)||(\xi_{h,n},Z_{h,n})||_\infty,
$$
(by definition of the predictor  step, using that $F$ is Lipschitz and $F(0)=0$), where $c$ 
is a constant depending only on  $\lambda$ and the Lipschitz constant of $F$. 
So, we have obtained our second important bound
\begin{equation}
\label{sup norm}
||(\xi_{h,n},Z_{h,n})||_\infty
\le ||(\xi_{0},Z_{0})||_\infty \exp(nhc),
\end{equation}
where $c$ 
depends only on $\lambda$ and the Lipschitz constant of $F$, and not on $\epsilon$.

\subsubsection*{Consistency estimates}

From now on, we fix
a time interval $[0,T]$ and denote by $c$ any constant depending only on the data
$T,X_0=I+\xi_0,Z_0,N,F,\epsilon,\lambda.$ 
Let $G$ be a smooth 1-periodic test function.
Let us fix a time step $n$. Since $G$ is 1-periodic, because of the rearrangement step, we have
\begin{equation}
\label{ineq}
\int_0^1 G(X_{h,n+1}(a))da
= \int_0^1 G(\hat X_{h,n+1}(a))da.
\end{equation}
Thus,  using the definition of scheme (\ref{TCM1}), we get
\begin{equation}
\label{key}
\begin{split}
\int_0^1 G(X_{h,n+1}(a))da = 
\\
\int_0^1 G(X_{h,n}(a)+h[\lambda \xi_{h,n}(a)+Z_{h,n}(a)]
+\sqrt{2\epsilon h}N(a/r))da.
\end{split}
\end{equation}
By Taylor expansion in $h$, 
we get for the right-hand side of (\ref{key}).
$$
\int_0^1 
\{G(X_{h,n}(a))
+G'(X_{h,n}(a))h(\lambda \xi_{h,n}(a)+Z_{h,n}(a))
+G'(X_{h,n}(a))\sqrt{2\epsilon h}N(a/r)
$$
$$
+\frac{1}{2}G"(X_{h,n}(a))[h(\lambda \xi_{h,n}(a)+Z_{h,n}(a))+\sqrt{2\epsilon h}N(a/r)]^2\}da+O(h^{3/2}),
$$
using the sup-norm estimate (\ref{sup norm}).
We want to simplify this expression, in order to obtain
$$
\int_0^1 
\{G(X_{h,n}(a))
+G'(X_{h,n}(a))h(\lambda \xi_{h,n}(a)+Z_{h,n}(a))
$$
$$
+G"(X_{h,n}(a))\epsilon h\}da+O(h^{3/2}),
$$
There are two terms to deal with:
$$
I_1=\int_0^1  
G'(X_{h,n}(a))\sqrt{2\epsilon h}N(a/r)da
$$
and
$$
I_2=\frac{1}{2}\int_0^1  
G"(X_{h,n}(a))
[h(\lambda \xi_{h,n}(a)+Z_{h,n}(a))+\sqrt{2\epsilon h}N(a/r)]^2 da,
$$
which both involve $\sqrt{h}$.
Let us write $N(a)=M'(a)$ where $M$ is a Lipschitz 1-periodic function $Q$,
so that $N(a/r)=r\partial_a[M(a/r)]$.
Thus,  integrating $I_1$ by part in $a$,
we get $I_1=-I_3+I_4$
 $$
I_3=-\int_0^1 
G"(X_{h,n}(a))\partial_a X_{h,n}(a)
\sqrt{2\epsilon h}rM(a/r)da
 $$
and the boundary term
$$
I_4=r \sqrt{2\epsilon h}[G'(X_{h,n}(1))M(1/r)-G'(X_{h,n}(0))M(0)].
$$
Since $1/r$ is assumed to be an integer, by 1-periodicity of $G$, $X_{h,n}-I$ and $M$, we first find $I_3=0$.
Then, using that $X_{h,n}(a)$ is non-decreasing in $a$, we get
$$
|I_3|\le cr\sqrt{h}\int_0^1 \partial_a X_{h,n}(a)da=cr\sqrt{h}(1+\int_0^1 \partial_a \xi_{h,n}(a)da),
$$
where $c$ denotes a generic constant depending on the data and the test function $G$.
Thanks to the sup norm estimate (\ref{sup norm}), we finally get
$|I_1|\le cr\sqrt{h}$.
Thus $I_1$ can be absorbed in the error term $O(h^{3/2})$, just by assuming 
$r=O(h)$.
Let us now consider $I_2=I_5+I_6+I_7$, where
$$
I_5=\frac{1}{2}\int_0^1  
G"(X_{h,n}(a))
[h(\lambda \xi_{h,n}(a)+Z_{h,n}(a))]^2 da
$$
$$
I_6=\int_0^1  
G"(X_{h,n}(a))
h(\lambda \xi_{h,n}(a)+Z_{h,n}(a))\sqrt{2\epsilon h}N(a/r)da
$$
$$
I_7=\int_0^1 \epsilon h G"(X_{h,n}(a))N^2(a/r)da
$$
Since we know by (\ref{sup norm}) 
that $(\xi_h,Z_h)$ is bounded in $L^\infty$, we get $I_5+I_6=O(h^2)$.
The last term, $I_7$ is just
$$
I_8=\int_0^1 \epsilon h G"(X_{h,n}(a))da,
$$
since $N^2(a)=1$.
Thus, assuming $r=O(h)$,
we have obtained the key consistency property of the scheme
\begin{equation}
\label{consistency discrete}
\begin{split}
\frac{1}{h}
\int_0^1 [G(X_{h,n+1}(a))-G(X_{h,n}(a))]da
=
\\
\int_0^1 
\{G'(X_{h,n}(a))(\lambda \xi_{h,n}(a)+Z_{h,n}(a))+\epsilon G"(X_{h,n}(a))\}da+O(\sqrt{h}),
\end{split}
\end{equation}
for all smooth 1-periodic function $G$.
We can write this relation in terms of the probability measure $\rho_{h,n}(dx)$ and the 
real-valued measure $Q_{h,n}(dx)$ defined by
\begin{equation}
\label{def rho Q}
\begin{split}
\int_0^1 g(x)\rho_{h,n}(dx)=\int_0^1 g(X_{h,n}(a))da,
\\
\int_0^1 g(x)Q_{h,n}(dx)=\int_0^1 g(X_{h,n}(a))(\lambda \xi_{h,n}(a)+Z_{h,n}(a))da,
\;\;\;\forall g\in C^0(\ETTE).
\end{split}
\end{equation}
Notice that both are uniformly bounded in $(h,n)$ as Borel measures, since $X_{h,n}$ and
$X_{h,n}$ are uniformly bounded, as $nh\le T$, for all finite $T$.
We obtain, as a substitute for (\ref{consistency discrete}),
\begin{equation}
\label{consistency discrete2}
\begin{split}
\frac{\rho_{h,n+1}-\rho_{h,n}}{h}+\partial_x Q_{h,n}=\epsilon \partial^2_{xx}\rho_{h,n}
+O(\sqrt{h}),
\end{split}
\end{equation}
in the sense of distributions on $\ETTE$.

\subsubsection*{Compactness}

We first perform a linear interpolation in time of the measures $(\rho_{h,n},Q_{h,n})$
which generates a continuous piecewise linear function of time 
$t\rightarrow (\rho_h(t),Q_h(t))$ valued in the dual space $C^0(\ETTE)'$ ( i.e. the space of all bounded Borel measures on $\ETTE$, by Riesz identification theorem).
For each fixed $T>0$, this family is bounded in 
$Lip([0,T],C^2(\ETTE)')$, because of (\ref{consistency discrete2}), 
and therefore compact in $C^0([0,T],C^0(\ETTE)')$.
Let us consider a convergent subsequence, for a sequence of time steps $h_m$
going to zero, and denote the limit $(\rho,Q)$.

Since $Z_h$ is "slaved" by $\xi_h$ and depends "linearly" on it through (\ref{TCM1}), we are just
concerned by the compactness of $\xi_h$.
Because of the sup norm and  Lipschitz estimates (\ref{sup norm},\ref{lipschitz}),
we already know that $(\xi_h(t))$ 
is valued in a fixed compact set of $L^2(\ETTE)$.
Next, we use in a crucial way the following lemma which comes from rearrangement theory
(or, alternately, from "optimal transport theory")  \cite{Br3,Vi}:
\begin{Proposition}
\label{compactness}
Let $(\xi_h)$ a family of functions in $H=L^2(\ETTE)$ such that
$X_h=I+\xi_h$ is a non-decreasing function.
Then $(\xi_h)$ converges to $\xi$ in $H$ if and only if
\begin{equation}
\label{convergencce}
\int_0^1 G(X_h(a))da
\rightarrow \int_0^1 G(X(a))da,
\end{equation}
where $X=I+\xi$, for all smooth 1-periodic unction $G$.
\end{Proposition}
Using this proposition, we deduce from (\ref{consistency discrete}) that 
$(\xi_h)$ is  uniformly equi-continuous from $[0,T]$ to $H$, for all $0<T<+\infty$.
We conclude, by the Arzela-Ascoli theorem
that the family $(\xi_h,Z_h)$ is relatively compact in $C^0(\ERRE_+,H)$, remembering that $Z_h$
is linearly slaved by $\xi_h$.

\subsubsection*{Consistency}

We now have enough compactness to pass to the limit: up to the extraction of a subsequence of
$h\rightarrow 0$, we have at least a limit $(\xi,Z)$ in $C^0(\ERRE_�,H)$ for
the family $(\xi_h,Z_h)$. Because of the sup-norm and Lipschitz estimates (\ref{sup norm},\ref{lipschitz}),
we also know that $(\xi(t,a),Z(t,a))$ are bounded and Lipschitz continuous in $a$, uniformly in $t$ for all bounded time interval $t\in [0,T]$.
Concerning $Z_h$, passing to the limit in the second equation of (\ref{TCM1}) is straightforward
(using that $F$ is smooth and Lipschitz). We get:
\begin{equation}
\label{Z evolution}
\frac{dZ}{dt}+\lambda Z=F(\xi)-\lambda^2 \xi.
\end{equation} 
Next, from (\ref{consistency discrete}), we deduce that
for all smooth 1-periodic function $G$ 
\begin{equation}
\label{consistency}
\begin{split}
\frac{d}{dt}\int_0^1 G(a+\xi(t,a))da=
\\
\int_0^1 \{G'(a+\xi(t,a))(\lambda \xi(t,a)+Z(t,a))+\epsilon G"(a+\xi(t,a))\}da.
\end{split}
\end{equation}

\subsubsection*{Uniqueness of the limit}

Let us consider any limit $(\xi,Z)$ of the scheme, on a finite interval of time $[0,T]$.
We know that it belongs to $L^\infty([0,T],Lip(\ETTE))$. Thus, for each $t$, $X(t,a)=a+\xi(t,a)$ has an inverse function $u(t,x)$ which is increasing and satisfies $u(t,x+1)=1+u(t,x)$. In addition, there is
a constant $\alpha>0$ such that  $\partial_x u(t,x)\ge \alpha$. (However, at this stage, $u(t,x)$ may have jumps in $x$). Its derivative $\rho(t,x)=\partial_x u(t,x)\ge\alpha$ is, for each $t$, a probability measure on $\ETTE$, bounded away from zero. 
We may also introduce, for each time $t$, the real valued bounded (uniformly in $t$) measure $Q(t,x)$ on $\ETTE$,
defined by
\begin{equation}
\label{Q}
\begin{split}
\int_0^1 g(x)Q(t,dx)=\int_0^1 g(X(t,a))(\lambda \xi(t,a)+Z(t,a))da,\;\;\;\forall g\in C^0(\ETTE),
\end{split}
\end{equation}
which is absolutely continuous with respect to $\rho(t,x)$ with a bounded density (since $X$ and $Z$ are bounded).
Thus, $(\rho,Q)$ belongs to $L^\infty([0,T],C^0(\ETTE)')$, where $C^0(\ETTE)'$ is the space of all bounded (Borel) measures on $\ETTE$.
The consistency relation (\ref{consistency}), written in terms of $\rho$ 
and $Q$, exactly means
\begin{equation}
\label{rho}
\begin{split}
\partial_t\rho+\partial_x Q=\epsilon\partial_{xx}\rho,
\end{split}
\end{equation}
in the sense of distribution. Their Fourier coefficients $F\rho(t,k)$, $FQ(t,k)$ are uniformly bounded in 
$t\in [0,T]$ and $k\in \ZED$, and satisfy
\begin{equation}
\label{rho fou}
\begin{split}
(\partial_t+4\pi^2\epsilon  k^2)F\rho(t,k)=2i\pi k FQ(t,k).
\end{split}
\end{equation}
This implies that $\rho$ belongs to $L^\infty([0,T],H^s(\ETTE))$ (where $H^s(\ETTE)$ denotes the standard Sobolev space of 1-periodic function with derivatives in $L^2$ up to order $s$, for $s\in \ERRE$) for  all $s<1/2$, and, therefore, by Sobolev embedding theorem to all
$L^\infty([0,T],L^q(\ETTE))$ for $1\le q<+\infty$. Since the density of $Q$ with respect to $\rho$ is bounded, we deduce that $Q\in L^\infty([0,T],L^q(\ETTE))$ for all $1\le q<+\infty$. 
Since $\rho=\partial_x u$, we also deduce $u\in L^\infty([0,T],H^s(\ETTE))$ for all $s<3/2$, and, thus
$u\in L^\infty([0,T],C^{0,\sigma}(\ETTE))$ for all $\sigma<1$. Thus, from the definition of $Q$ (\ref{Q}), we can now write
\begin{equation}
\label{Qbis}
\begin{split}
\rho(t,x)=\partial_x u(t,x),\;\;
Q(t,x)=((x-u(t,x))\lambda +Z(t,u(t,x))\partial_x u(t,x)
\end{split}
\end{equation}
and deduce that $Q$, just like $\rho$, belongs to $L^\infty([0,T],H^s(\ETTE))$ for all $s<1/2$. 
[For that, instead of using Fourier coefficients, we estimate translations of $Q$ in $L^2$. Using that 
$u$ is bounded and $Z$ Lipschitz, we get
$$
I=\int_0^1 |Q(t,x+\omega)-Q(x)|^2 dx\le (I_1+I_2)c
$$
where 
$$
I_1=\int_0^1 |\rho(t,x+\omega)-\rho(x)|^2 dx
$$
$$
I_2=\sqrt{\int_0^1 |u(t,x+\omega)-u(x)|^4 dx \int_0^1 \rho(t,x)^4 dx}.
$$
Since $\rho$ belongs to $L^\infty([0,T],H^s(\ETTE))$ for all $s<1/2$, we immediately get that
\\
$I_1\le c_s |\omega|^{2s}$, for all $s<1/2$. Next, we use that $u$ belongs to $L^\infty([0,T],H^s(\ETTE))$ for all $s<3/2$. This implies by Sobolev embedding theorem, that $u$ and $Z(u)$ certainly belong to 
$L^\infty([0,T],W^{1,4}(\ETTE))$, while $\rho$ belongs to $L^\infty([0,T],L^4(\ETTE))$. 
This implies that $I_2\le c|\omega|^2$. So we see that $I\le c_s|\omega|^{2s}$ for all $s<1/2$,
which, indeed, means that $Q$ belongs to $L^\infty([0,T],H^s(\ETTE))$ for all $s<1/2$.]
So we can now bootstrap the regularity of $\rho$ and $Q$, using (\ref{rho fou}).
We first get that $\rho$, belongs to $L^\infty([0,T],H^s(\ETTE))$ for all $s<3/2$ and $u$ to
the this space for all $s<5/2$. By Sobolev embedding theorem, this means that
$u$ also belong to $L^\infty([0,T],C^{1,\sigma}(\ETTE))$ for all $\sigma<1$. We already know that
$\partial_x u(t,x)\ge \alpha>0$. 
Since $X(t,a)$ is the inverse function of  $u$, it follows that 
$\xi(t,a)=X(t,a)-a$ also belongs to
$L^\infty([0,T],C^{1,\sigma}(\ETTE))$ for all $\sigma<1$. This is also true for $Z$, which is slaved by 
$\xi$ through (\ref{Z evolution}).
 At this point, the boostrap argument becomes pure routine and we conclude that
$u,\xi,Z$ are smooth, just as the initial conditions.
From (\ref{rho},\ref{Qbis},\ref{Z evolution}), we see they are the unique
solution in classical sense of the system
\begin{equation}
\begin{split}
u(t,X(t,a))=a,\;\;\;
X(t,a)=a+\xi(t,a),
\\
\partial_t u(t,x)+[Z(t,u(t,x))+\lambda (x-u(t,x))]\partial_x u(t,x)=\epsilon\partial^2_{xx} u,
\\
\partial_t Z(t,a)+\lambda Z(t,a)=F(\xi(t,a))-\lambda^2 \xi(t,a).
\end{split}
\end{equation}
This concludes the proof of our theorem since this system is
just the formulation (\ref{scalar3}) of the Navier-Stokes system
(\ref{NSP material}).

\subsection*{Acknowledgment}
This work is partly supported by the ANR grant OTARIE
ANR-07-BLAN-0235. The author thanks the Dipartimento di Matematica Pura e
Applicata of the Universit\`a degli Studi dell'Aquila, l'Aquila, Italy, where this work was partly done.

\end{document}